\documentstyle{amsppt}
\tolerance 3000
\pagewidth{5.5in}
\vsize7.0in
\magnification=\magstep1
\widestnumber \key{AAAAAAA} 
\topmatter
\author Alex Iosevich and Steen Pedersen     
\endauthor 
\thanks Research supported in part by the NSF grant DMS97-06825 
\endthanks 
\address Department of Mathematics Wright State University Dayton Ohio 
45435 e-mail: iosevich$\@$zara.math.wright.edu steen$\@$math.wright.edu
\endaddress 
\abstract Let $D$ be a bounded domain in ${\Bbb R}^n$ whose boundary has 
a Minkowski dimension $\alpha<n$. Suppose that $E_{\Lambda}=
{\{e^{2 \pi i x \cdot \lambda}\}}_{\lambda \in \Lambda}$, $\Lambda$ an 
infinite discrete subset of ${\Bbb R}^n$, is a frame of exponentials for
$L^2(D)$, with frame constants $A,B$, $A \leq B$. Then if 
$$ R \ge C{\left(\frac{{B|\partial D|}_{\alpha}}{A|D|} \right)}^
{\frac{1}{n-\alpha}},$$
where $C$ depends only on the ambient dimension $n$ and 
${|\partial D|}_{\alpha}$ denotes the Minkowski content, then every cube of 
sidelength $R$ contains at least one element of $\Lambda$. We give examples
that illustrate the extent to which our estimates are sharp. \endabstract 
\title How large are the spectral gaps?  
\endtitle
\endtopmatter 

Let $D$ be a domain of finite Lebesgue measure in ${\Bbb R}^n$. 
Suppose that $L^2(D)$ has a frame of exponentials of the form 
$E_{\Lambda}=\{e^{2 \pi i x \cdot \lambda}\}$, $\lambda \in \Lambda$, a 
discrete infinite subset of ${\Bbb R}^n$, with frame constants $A,B$, 
$A \leq B$, in the sense that 
$$ A{||f||}_{L^2(D)}^2 \leq \sum_{\Lambda} {|\hat{f}(\lambda)|}^2 
\leq B{||f||}_{L^2(D)}^2, \tag*$$ where $f:D \rightarrow {\Bbb C}$, and 
$\hat{f}$ denotes the standard Fourier transform. In this paper we will 
work with frames rather than exponential basis because $L^2$ of every 
bounded domain has frames, whereas exponential basis are generally hard to 
come by. (See \cite{Fug}). The following quantities were introduced by 
Beurling. See \cite{Br}. 
 
$$ D_R^{+}= \max_{x \in {\Bbb R}^n} \# \{ \Lambda \cap Q_R(x)\}, \tag1$$ 
where $Q_R(x)$ is a cube of sidelength $2R$ centered at $x$, and let 
$$ D_R^{-}= \min_{x \in {\Bbb R}^n}\#\{ \Lambda \cap Q_R(x)\}. \tag2$$ 

It follows from results proved by Landau (\cite{Lan}, see also 
\cite{GR}) that if $D$ is a bounded domain then  
$$ \limsup_{R \rightarrow \infty} \frac{D_R^{-}}{{(2R)}^n} \ge |D|. \tag3$$ 

If the set $E_{\Lambda}$ is actually an orthogonal basis for $L^2(D)$ then 
the inequality $(3)$ is actually an equality for both $D_R^{-}$ and $D_R^{+}$.  

These results show that, asymptotically, a sufficiently large cube centered
at any point contains the number of elements of $\Lambda$ at least equal 
to its volume multiplied by the Lebesgue measure of $D$. In this paper we 
will show that if the Minkowski dimension, $\alpha$,  of the boundary 
$\partial D$ is smaller than the ambient dimension $n$, then there exists 
$$R=C{\left(\frac{{B|\partial D|}_{\alpha}}{A|D|}\right)}^{\frac{1}{n-\alpha}},
\tag4$$ where $C$ only depends on $n$ and
${|\partial D|}_{\alpha}=\overline{\lim}_{\epsilon \downarrow 0} 
{\epsilon}^{\alpha-n} |\{x: d(x, \partial D)<\epsilon\}|$ denotes the 
$\alpha$-dimensional upper Minkowski content of $\partial D$, 
such that a cube of sidelength $2R$ centered at any point contains at least
one element of $\Lambda$. Note that if $\partial D$ is, say, piecewise smooth,
then $\alpha=n-1$ and $R=C\frac{B|\partial D|}{A|D|}$. 

A note on notation. The letter $C$ below shall denote a generic constant 
which may change from line to line. We shall give more precise information
about the constants when appropriate. 

Our main result is the following. 

\proclaim{Theorem 1} Let $D$ denote a domain in ${\Bbb R}^n$ with finite 
non-zero Lebesgue measure whose boundary $\partial D$ has Minkowski dimension 
$\alpha<n$ in the sense that 
$$ |\{x \in {\Bbb R}^n: d(x, \partial D)< \epsilon\}| \leq 
C {\epsilon}^{n-\alpha}. \tag5$$

Then there exists $C$ depending only on $n$, such that if 
$$ R \ge C 
{\left(\frac{{B|\partial D|}_{\alpha}}{A|D|}\right)}^{\frac{1}{n-\alpha}},
\tag6$$ then  
$$ \Lambda \cap Q_R(\mu) \not= \emptyset \tag7$$ for every  
$\mu \in {\Bbb R}^n$, and any set $\Lambda$ such that $E_{\Lambda}$ is an 
exponential frame for $L^2(D)$, with frame constants $A,B$, $A\leq B$ 
where $Q_R(\mu)$ denotes the cube of sidelength $2R$ centered at $\mu$. 
\endproclaim 

In other words, our result shows, at least if $A=B$, 
that if $D$ has a fixed volume, then the 
maximum gap between the elements of $\Lambda$ is no larger than a fixed
constant times the the Minkowski content of the boundary. Moreover, 
the rate of increase depends on the Minkowski dimension of $\partial D$. 
This idea is illustrated by the following simple example. 

\proclaim{Example 2} Let $D=[0,a_1] \times [0,a_2] \times \dots \times 
[0,a_n]$, $a_1 \ge a_2 \ge \dots \ge a_n>0$, $\Pi_{j=1}^n a_j=1$. 
We can take $\Lambda=\Pi_{j=1}^n 
\frac{1}{a_j} {\Bbb Z}$. It is not hard to see that the largest cube that
does not intersect $\Lambda$ has sidelength $2R=\frac{1}{a_n}$. The measure
of $\partial D$ is $2\sum_{j=1}^n \frac{1}{a_j}$. It follows that 
$$ \frac{1}{4n} \leq \frac{R}{|\partial D|} \leq \frac{1}{4}, \tag8$$ 
so $R$ grows linearly with $|\partial D|$. 
\endproclaim 

\proclaim{Example 3} We now spice up the above example to illustrate the 
fractal phenomenon. Let $D$ be a domain constructed by taking a square 
${[0,1]}^2$ and replacing the upper and lower segments  by identical 
fractal curves of Minkowski dimension $1<\alpha<2$. It is not hard to see
that $\Lambda$ may be taken to be ${\Bbb Z}^2$. (See \cite{Fug}). We
now blow up the domain by the factor of $t>1$ (i.e we apply the matrix $tI$, 
where $I$ is the identity matrix). Let $tD$ denote the image of $D$ under 
that mapping. The set $\Lambda$ must now be taken to be 
${(\frac{1}{t}\Bbb Z)}^2$, which tells us that $R$ in Theorem 1 should be 
$\approx \frac{1}{t}$. On the other hand, ${|\partial tD|}_{\alpha}
\approx t^{\alpha}$, and $|tD|=t^2$, so Theorem 1 gives us 
$R \approx {\left( \frac{t^{\alpha}}{t^2} \right)}^{\frac{1}{2-\alpha}}=
\frac{1}{t}$. \endproclaim 

The following example shows that if the Lebesgue measure $|D|=0$ the 
conclusion of Theorem 1 no longer holds. 
\proclaim{Example 4} Let $D \subset [0,1]$ denote the Cantor type set 
consisting of numbers that do not have $1$ or $3$ in their base $4$ 
expansion. Let $m$ denote the unique probability measure supported on $D$ 
(see \cite{Fal}) given by the equation 
$$ \int f(t) dm(t)=\frac{1}{2} \int f\left(\frac{t}{4}\right) dm(t)+
\frac{1}{2} \int f\left(\frac{t+2}{4}\right) dm(t). \tag9$$ 

One can check that
$$ \widehat{m}(t)=e^{\pi i \frac{2}{3}t} \Pi_{j=0}^{\infty} 
\cos\left(\frac{\pi t}{2 \cdot 4^n} \right). \tag10$$ 

If $\Lambda$ is the set of non-negative integers whose 
base $4$ expansion does not contain $2$ or $3$, then $E_{\Lambda}$ is an 
orthonormal basis of $L^2(m)$. (See \cite{JP}).   

In particular this shows that the conclusion of Theorem 1 fails miserably 
in this case.
\endproclaim  

\proclaim{Example 5} In this example we shall see that there exist families
of domains with piecewise smooth boundaries such that the volume of each 
domain is $1$, the length of the boundary tends to infinity, but $R$, in 
the sense of Theorem 1, may always be taken to be $\frac{1}{2}+\epsilon$, 
for any $\epsilon>0$. 

Let $D_k$ denote the unit square in ${\Bbb R}^2$ where the upper and lower
edges are replaced by a sawtooth function with $k$ teeth where the height
of each tooth is $\frac{1}{2}$. The length $|\partial D_k|$ goes to 
infinity as $k \rightarrow \infty$. The set $\Lambda$ for each $D_k$ is 
${\Bbb Z}^2$, so $R$, in the sense of Theorem 1, may always be taken to be 
$\frac{1}{2}+\epsilon$, for any $\epsilon>0$. This says that the inequality
$(6)$ does not sharply describe the behavior of $R$ in this case. However, 
the proof of Theorem 1 (see the discussion at the end of the proof of Theorem 
1 below) shows that in some cases $R$ may be taken to be 
$C \frac{diameter(D)}{|D|}$, where $C$ depends only on $n$. We shall see that
the example given in this paragraph falls into that category. 
\endproclaim  

In all the previous examples we used frames which were actually orthogonal 
bases. However, this phenomenon persists in the cases when orthogonal 
exponential basis do not exist and we have to make do with frames. 

\proclaim{Example 6} Let $B_r$ denote the disc of radius $r$ in ${\Bbb R}^2$
centered at the origin. It was shown in \cite{JP2} that 
$\Lambda=\frac{1}{2r} {\Bbb Z}^2$ is frame for $L^2(B_r)$ with constants 
$A=B=4r^2$. Note that we do not have orthogonal basis becuase, in particular,
that would imply that $A=B=|B_r|=\pi r^2$. It is well known that $B_r$ 
does not have orthogonal basis of exponentials. See \cite{Fug}. 

It is clear that $R$, in the sense of Theorem 1 must be taken to be 
greater than $\frac{1}{4r}$, which is exactly what Theorem 1 predicts. 
\endproclaim 

The key estimate (see Lemma 9 below) involved in the proof of Theorem 1 is
$$ \sum_{\lambda \notin Q_R(\mu)} {|\widehat{\chi}_D(\lambda-\mu)|}^2 \leq 
C\frac{{|\partial D|}_{\alpha}}{R^{n-\alpha}}, \tag11$$ for any 
$\mu \in {\Bbb R}^n$, where $C$ depends only on the dimension and 
on the frame constant $B$. 

This estimate is similar to the estimate that comes up in the theory of 
irregularities of distributions, (see \cite{Mgr}, p.110),  namely that for any 
domain $S$ whose boundary is a piecewise $C^1$ curve ${\Cal C}$ 
$$ \int_{|t| \ge R} {|\widehat{\chi}_S(t)|}^2 dt \leq \frac{|{\Cal C}|}
{2 \pi^2 R}. \tag12$$ 

In fact, our proof of the estimate $(11)$ given in Lemma 9 below uses an idea 
from the proof of the estimate $(12)$ given by Brandolini, Colzani, and 
Travaglini in \cite{BCT}.  

The proof of Theorem 1 is based on the following sequence of lemmae. 

\proclaim{Lemma 7} For any $f \in L^2(D)$ define 
$$ F_Df(\xi)=\int_D e^{-2\pi ix \cdot \xi} f(x)dx, \tag13$$ and let
$\hat{f}$ denote the standard Fourier transform 
$$ \hat{f}(\xi)=\int_{{\Bbb R}^n} e^{-2\pi i x \cdot \xi} f(x)dx. \tag14$$ 

Let $t_hf(x)=f(x+h)$, and let $\chi_D$ denote the characteristic function of 
$D$. Then 
$$ F_Dt_h \chi_{D}(\lambda) =e^{2\pi i \lambda \cdot h} 
\widehat{\chi}_{D \cap D+h}(\lambda), \tag15$$ 
$$ F_D t_{-h}\chi_D(\lambda)=\widehat{\chi}_{D \cap D+h}(\lambda), \tag16$$ and
$$ F_D\chi_D(\lambda)=\widehat{\chi}_D(\lambda). \tag17$$ 
\endproclaim 

The proof is straightforward. 

\proclaim{Lemma 8} Let $D$ be as above. Then 
$$ \int_{D} {|\chi_D(x+h)-\chi_D(x-h)|}^2 dx \leq C{|h|}^{n-\alpha}, 
\tag18$$ and
$$ \int_{D} {|\chi_D(x)-\chi_D(x-h)|}^2 dx \leq C{|h|}^{n-\alpha}, \tag19$$ 
with $C \leq C'{|\partial D|}_{\alpha}$, where $C'$ depends only on $n$.  
\endproclaim 

\remark{Remark} We note again that even though the estimate 
$C \leq C'{|\partial D|}_{\alpha}$ is best possible over all $h$'s, for 
special choices of $h$, the estimate is frequently much better. (See 
Example 5 above). \endremark 

To prove $(19)$ note that the left hand side equals 
$|\{D-(D+h)\} \cup \{(D+h)-D\}| \leq |\{x \in {\Bbb R}^n: 
d(x, \partial D)<h\}| \leq C {|\partial D|}_{\alpha} {|h|}^{n-\alpha}$. 
The proof of $(18)$ is similar. 

The key lemma is the following. (See \cite{BCT} for a similar argument). 
\proclaim{Lemma 9} Let $D$ be as above and let $\Lambda$ be such that 
$E_{\Lambda}$ is a frame of $L^2(D)$ with frame constants $A$ and $B$, 
$A \leq B$. Then 
$$ \sum_{\{ \lambda \in Q_{2^{k+1}}-Q_{2^k}\}} {|\widehat{\chi}_D(\lambda)|}^2 
\leq CB2^{-k(n-\alpha)}, \tag20$$ where $Q_R=Q_R(0, ..., 0)$, and $C$ as in 
Lemma 8. 
\endproclaim 

To prove Lemma 9 chose $N$ boxes $A_k^j$ and $N$ vectors $h_j$ such that 
$2^{-k} \leq |h_j| \leq 2^{-k+1}$, $\cup_{j=1}^N A_k^j=Q_{2^{k+1}}-Q_{2^k}$,
and 
$$ \left|e^{2 \pi i \lambda \cdot h_j}-1 \right| \ge 1, \ \ 
\lambda \in A_k^j. \tag21$$ 

Clearly this can be done in any dimension $n$, for a sufficiently large 
$N=N(n)$. 

Now, by triangle inequality 
$$ {\left( \sum_{A_k^j} {|\widehat{\chi}_D(\lambda)|}^2 \right)}^{\frac{1}{2}} 
\leq {\left( 
\sum_{A_k^j} {|\widehat{\chi}_{D \cap D+h_j}(\lambda)|}^2 \right)}^
{\frac{1}{2}}+ {\left(
\sum_{A_k^j} {|\widehat{\chi}_D(\lambda)-\widehat{\chi}_{D \cap D+h_j} 
(\lambda)|}^2 \right)}^{\frac{1}{2}}=I+II. \tag22$$ 

By Lemma 7, the frame property, and Lemma 8 we get 
$$ {II}^2 \leq \sum_{\Lambda} {|\widehat{\chi}_D(\lambda)-\widehat{\chi}_
{D \cap D+h_j}(\lambda)|}^2=$$ 
$$ \sum_{\Lambda} {|F_D\chi_D(\lambda)-F_Dt_{-h_j}\chi_D(\lambda)}|^2 \leq$$ 
$$ B \int_D {|\chi_D(x)-\chi_D(x-h_j)|}^2 dx \leq 
CB{|h_j|}^{n-\alpha} \leq CB2^{-k(n-\alpha)}. \tag23$$

On the other hand, by $(21)$, Lemma 7, the frame property, and Lemma 8 we get
$$ I^2 \leq \sum_{A_k^j} {|\widehat{\chi}_{D \cap D+h_j}(\lambda)|}^2 
{|e^{2 \pi i \lambda \cdot h}-1|}^2 \leq 
\sum_{\Lambda} {|\widehat{\chi}_{D \cap D+h_j}(\lambda)|}^2 
{|e^{2 \pi i \lambda \cdot h_j}-1|}^2=$$ 
$$ \sum_{\Lambda} {|F_Dt_{h_j}\chi_D(\lambda)-F_Dt_{-h_j}\chi_D(\lambda)|}^2
\leq$$ 
$$ B \int_D {|\chi_D(x+h_j)-\chi_D(x-h_j)|}^2 dx \leq 
CB{|h_j|}^{n-\alpha} \leq CB2^{-k(n-\alpha)}. 
\tag24$$ 

\head Proof of Theorem 1 \endhead 
\vskip.125in 

Since $E_{\Lambda}$ is a frame for $L^2(D)$ if and only if 
$E_{\Lambda-\mu}$ is also a frame for $L^2(D)$ (with the same frame 
constants) for any 
$\mu \in {\Bbb R}^n$, and our estimates do not 
depend on the choice of $\Lambda$, it is sufficient to 
consider the case $\mu=(0, \dots, 0)$. 

By the frame property and Lemma 7 we get 
$$ A|D| \leq  \sum_{\Lambda} {|F_D\chi_D(\lambda)|}^2=
\sum_{\Lambda} {|\widehat{\chi}_D(\lambda)|}^2= 
\sum_{Q_R} {|\widehat{\chi}_D(\lambda)|}^2 + 
\sum_{\{\lambda \notin Q_R\}} {|\widehat{\chi}_D(\lambda)|}^2. \tag25$$ 

Using Lemma 9 we see that if $R=2^{k_0}$,  
$$ \sum_{\{\lambda \notin Q_R\}} {|\widehat{\chi}_D(\lambda)|}^2 =
\sum_{k=k_0}^{\infty} \sum_{Q_{2^{k+1}}-Q_{2^k}} 
{|\widehat{\chi}_D(\lambda)|}^2 \leq 
CB2^{-k_0(n-\alpha)}=\frac{BC}{R^{n-\alpha}}. \tag26$$ 

So by $(25)$ and $(26)$
$$ \sum_{Q_R} {|\widehat{\chi}_D(\lambda)|}^2 \ge 
A|D|-\frac{BC}{R^{n-\alpha}} \tag27$$ 

which proves that if $R>{\left(\frac{BC}{A|D|}\right)}^{\frac{1}{n-\alpha}}$, 
then 
$$ \Lambda \cap Q_R \not=\emptyset. \tag28$$ 

Moreover, the above proof shows that 
$C \leq C'{|\partial D|}_{\alpha}$ where $C'$ 
depends only on $n$.
\remark{Remark} In the proof above the key estimate is 
$|\{D \cap D-h\}| \leq C{|h|}^{n-\alpha} {|\partial D|}_{\alpha}$. While 
this is the best possible estimate uniform in $h$, in the proof we are   
have a wide choice of $h$'s as long as $|h|=2^{-k}$ and the estimates 
$(18)$, $(19)$, and $(21)$ are satisfied. 

This observation can be used to handle the family of examples given by 
Example 5 above. For convenience we take $\Lambda=(\frac{1}{2}, 0)+
{\Bbb Z}^2$. We can now take all $h$'s in the proof of Theorem 1 of the 
form $h=(h_1,0)$ and for this choice of $h$'s it is easy to check that
$|\{D_k \cap D_k-h\}| \leq C|h| diameter(D_k)$, where $C$ is a uniform constant,
since the "teeth" of $D_k$'s point in the $y$-direction. Since 
$diameter(D_k)$ is uniformly bounded above and below, the lack of  
sharpness of Theorem 1 exposed in Example 5 is resolved for this family 
of examples. \endremark 

\newpage 

\centerline{References} 
\vskip.125in 

\ref \key Br \by A. Beurling \paper Local harmonic analysis with some 
applications to differential operators \yr 1966 \jour Some recent advances 
in the basic sciences, Academic press \vol 1 \endref 

\ref \key BCT \by L. Brandolini, L. Colzani, and G. Travaglini \paper 
Average decay of Fourier transforms and integer points in polyhedra 
\jour Ark. Mat. \vol 35 \yr 1997 \pages 253-275 \endref 

\ref \key Fal \by K. Falconer \paper The geometry of fractal sets 
\jour Cambridge University Press \yr 1986 \endref 

\ref \key Fug \by B. Fuglede \paper Commuting self-adjoint partial 
differential operators and a group theoretic problem \jour J. Funct. Anal. 
\yr 1974 \vol 16 \pages 101-121 \endref 

\ref \key GR \by K. Gr\"{o}chenig and H. Razafinjatovo \paper on Landau's 
necessary density conditions for sampling and interpolation of band-limited
functions \jour  J. London. Math. Soc. \vol 54 \pages 557-565 \yr 1996 
\endref 

\ref \key IP \by A. Iosevich and S. Pedersen \paper Spectral and tiling 
properties of the unit cube \jour (submitted for publication) \yr 1998 
\endref 

\ref \key JP \by P. E. T. Jorgensen and S. Pedersen \paper Dense analytic
subspaces in fractal $L^2$-spaces \jour J. Anal. Math. (to appear) 
\yr 1998 \endref 

\ref \key JP2 \by P. E. T. Jorgensen and S. Pedersen \paper Local
harmonic analysis for domains in $\Bbb{R}^n$ of finite measure
\jour Analysis and Topology, Eds. C. Andreian Cazacu, O. Lehto and
Th. M. Rassias, World Scientific Publishing Co. \yr 1997 \endref 

\ref \key LRW \by J. Lagarias, J. Reed, and Y. Wang \paper Orthonormal 
bases of exponentials for the $n$-cube \jour (preprint) \yr 1998 \endref 

\ref \key Lan \by H. Landau \paper Necessary density conditions for 
sampling and interpolation of certain entire functions \jour Acta Math. 
\vol 117 \pages 37-52 \yr 1967 \endref 

\ref \key Mgr \by H. Montgomery \paper Ten lectures on the interface between
analytic number theory and harmonic analysis \jour CBMS Regional conference 
series in mathematics \yr 1994 \endref 
\enddocument